\documentclass[psamsfonts]{amsart}

\usepackage{amssymb,amsfonts}
\usepackage[all,arc]{xy}
\usepackage{enumerate}
\usepackage{mathrsfs}

\theoremstyle{definition}

\theoremstyle{remark}

\makeatletter
\let\c@equation\c@thm
\makeatother
\numberwithin{equation}{section}

\title{On Cartesian Products which Determine Few Distinct Distances}

\author{Cosmin Pohoata}

\begin{document}

\begin{abstract}
Every set of points $\mathcal{P}$ determines $\Omega(|\mathcal{P}| / \log |\mathcal{P}|)$ distances. A close version of this was initially conjectured by Erd\H{o}s in 1946 and rather recently proved by Guth and Katz. We show that when near this lower bound, a point set $\mathcal{P}$ of the form $A \times A$ must satisfy $|A - A| \ll |A|^{2-\frac{2}{7}} \log^{\frac{1}{7}} |A|$. This improves recent results of Hanson and Roche-Newton. 
\end{abstract}

\maketitle

\section{Introduction}

\bigskip

Let $\mathcal{P}$ be a set of points in plane, and let $\Delta(\mathcal{P})$ denote the set of squares of distances spanned by $\mathcal{P}$. In other words,
$$\Delta(\mathcal{P}) = \left\{(p_{1}-q_{1})^{2} + (p_{2}-q_{2})^{2}: (p_{1},p_{2}), (q_{1},q_{2}) \in \mathcal{P}\right\}.$$
In [7], Guth and Katz showed that $\Delta(\mathcal{P}) \gg |\mathcal{P}| / \log |\mathcal{P}|$, where $\gg$ represents the usual Vinogradov symbol. When $P = A \times B$ for some finite sets of reals $A$ and $B$,  $\Delta(A \times B) = (A-B)^2 + (A-B)^2$, so this says that
$$|(A-B)^2 + (A-B)^2| \gg \frac{|A||B|}{\log |A||B|}.$$
In [5], Erd\H{o}s originally conjectured that all sets $\mathcal{P}$ should determine $\Omega(|\mathcal{P}| / \sqrt{\log |\mathcal{P}|})$ distinct distances, so the Guth-Katz bound is almost optimal. Nonetheless, very little is known for sets that achieve this bound. It is widely believed that sets with $O(|\mathcal{P}| / \log |\mathcal{P}|)$ distinct distances should come from some type of lattice. This is very well-motivated by the following beautiful result of Bernays [2], which generalizes a classical theorem of Landau.

\bigskip

{\bf{Theorem 1}}. {\it{Let $f(x,y) = ax^{2}+bxy+cy^{2}$ for integers $a,b,c \in \mathbb{Z}$, such that the determinant $b^{2}-4ac$ is not an integer square. Then, the number of integers between $1$ and $n$ that can be expressed as $f(u,v)$ with $u,v \in \mathbb{Z}$ is $O(n / \sqrt{\log n})$.}}

\bigskip

Using Theorem 1, one can easily check that sets with $O(n / \sqrt{\log n})$ distinct distances are given by $\sqrt{n} \times \sqrt{n}$ subsets of the integer lattice, the (equilateral) triangular lattice, or, more exotically, by the rectangular lattice
$$\mathcal{L}_{r} = \left\{(i,j\sqrt{r})\ |\ i,j \in \mathbb{Z},\ 1 \leq i, j \leq n\right\},$$
for every integer $r>1$. We refer the reader to [14] for a more detailed presentation of this discussion, where Sheffer also points out that unlike the first two examples, the latices $\mathcal{L}_{r}$ do not span squares or equilateral triangles.

In this paper, we will only take a look at sets that come from cartesian products, and show that whenever they determine few distinct distances they must exhibit some additive structure. Specifically, when $\mathcal{P}=A \times A$, we show that when the Guth-Katz bound is close to being tight, we have that
$$|A - A| \ll |A|^{2-\frac{2}{7}} \log^{\frac{1}{7}} |A|.$$
In light of the bipartite distance problems discussed by Brunner and Sharir in [3] and by Sheffer and the author in [12], we also consider the problem of showing that if there are few distinct distances between two cartesian products $A \times A$ and $B \times B$, then one of $A$ or $B$ has additive structure. We state both of these results more formally below.

\bigskip

{\bf{Theorem 2}}. {\it{Suppose $A$ is a finite set of real numbers and let $\Delta(A \times A)$ be the set of distances spanned by $A \times A$. Then,
$$|A-A| \ll |\Delta(A \times A)|^{\frac{6}{7}} \log^{\frac{1}{7}} |A|,$$
or equivalently $|D| \ll |D^{2}+D^{2}|^{\frac{6}{7}} \log^{\frac{1}{7}} |D|$, where $D$ denotes the difference set $A-A$.}}

\bigskip

{\bf{Theorem 3}}.  {\it{Suppose $A$ and $B$ are finite sets of real numbers and let $\Delta(A \times A, B \times B)$ be the set of distances between points in $A \times A$ and points in $B \times B$. Then,
$$\min\left\{|A-A|, |B-B|, |A-B|\right\} \ll |\Delta(A \times A, B \times B)|^{1-\frac{13}{205}} \cdot \L(A,B),$$
where
$$\L(A,B)=\min\left\{ \log^{\frac{3}{205}} |A|, \log^{\frac{3}{205}} |B|\right\}.$$}}

\bigskip

In particular, if $|\Delta(A \times A)| \ll |A|^{2}$ holds in Theorem 2, then $|A-A| << |A|^{2-\frac{2}{7}} \log^{\frac{1}{7}} |A|$. The improves a recent theorem by Hanson [8], who showed that under this hypothesis we have that $|A-A| \ll |A|^{2-\frac{1}{8}}$. In the meantime this was also sharpened by Roche-Newton in [13], who showed $|A-A| \ll |A|^{2-\frac{2}{11}}$, but Theorem 2 is stronger. Our proof will rely on the sum-product estimate of Solymosi from [15] as a black-box:

\bigskip

{\bf{Theorem 5}}. {\it{Let $S \subset \mathbb{R}$ be a set. Then,
$$|S+S|^{2}|SS| \geq \frac{|S|^{4}}{4\lceil \log |S|\rceil}.$$}}

\bigskip

The proof of Theorem 3 will rely on two results. The first one is the following Lemma by Balog [1], which comes from Solymosi's original idea for Theorem 5.

\bigskip

{\bf{Lemma 6}}. {\it{Let $R$, $S$, $T$ be finite sets of real numbers. Then
$$|RT+RT| |ST+ST| \gg |R/S||T|^2.$$}}

\bigskip

The second one is the following Lemma due to Shkredov, which is Theorem 3 in [15] (and the statement of which should be in some sense compared to that of Theorem 2 above).

\bigskip

{\bf{Lemma 7}}. {\it{Let $A \subset \mathbb{R}$ be a finite set and let $D = A-A$. Then
$$|D/D| \gg |D|^{1+\frac{1}{12}} \log^{-1/4}|D|.$$}}

\bigskip

Last but not least, we will also need the classical Pl$\ddot{\text{u}}$nnecke-Ruzsa inequality, for which a simple proof can be found in [11].

\bigskip

{\bf{Lemma 8}}. {\it{Let $A \subset \mathbb{R}$ be a finite set. Then
$$|k \cdot A - \ell \cdot A| \leq \frac{|A+A|^{k+\ell}}{|A|^{k + \ell - 1}}.$$}}

\bigskip

\section{Proof of Theorem 2}

\bigskip

If $D = A-A$, then $|\Delta(A \times A)| = |D^{2}+D^{2}|$, where $D^{2} = \left\{(x-y)^2: x,y \in A\right\}$. We claim that
$$|D| \ll |D^{2}+D^{2}|^{\frac{6}{7}} \log^{\frac{1}{7}} |D|.$$
We apply Theorem 2 for the set $S:=D^{2}$. Using the observation that $|D^{2} D^{2}|$ is equal to $|DD|$ (up to a small constant), this yields
$$|D^2+D^2|^{2}|DD| \geq  |D^2+D^2|^{2}|D^{2} D^{2}| \geq \frac{|D^{2}|^{4}}{4\lceil \log |D^{2}|\rceil} \gg \frac{|D|^{4}}{\log |D|}.$$
On the other hand for every four real numbers $a_{1},a_{2},b_{1},b_{2}$, we have that
$$(b_{1}-a_{1})^2 + (b_{2}-a_{2})^2 - (b_{1}-a_{2})^2 - (b_{2}-a_{1})^2 = 2(a_{2}-a_{1})(b_{1}-b_{2}),$$
which yields the inclusion
$$2DD \subset 2\cdot D^2 - 2\cdot D^2.$$ 
We emphasize here that for $X \subset \mathbb{R}$ and $c \in \mathbb{Z}_{>0}$, the set $cX$ denotes the set of scalar multiples $ \left\{cx: x \in X\right\}$, whereas $c \cdot X$ denotes the sumset $\sum_{i=1}^{c} X$. The inclusion together with Lemma 8 then yield
\begin{eqnarray*}
|D^2+D^2|^{2}|DD| &=& |D^2+D^2|^{2}|2DD|\\
&\leq& |D^{2}+D^{2}|^{2} |2 \cdot D^{2}-2 \cdot D^{2}|\\
&\ll& |D^{2}+D^{2}|^{2} \left( \frac{|D^{2}+D^{2}|^{4}}{|D|^{3}} \right).
\end{eqnarray*}
Putting the two bounds together, we conclude that
$$\frac{|D^{2}+D^{2}|^{6}}{|D|^{3}} \gg \frac{|D|^{4}}{\log |D|},$$
which yields
$$|D| \ll |D^{2}+D^{2}|^{\frac{6}{7}} \log^{\frac{1}{7}} |D|.$$
$\hfill \square$

\bigskip

\section{Proof of Theorem 3}

\bigskip

For convenience, write again that $|\Delta(A \times A, B \times B)|=|(A-B)^2+(A-B)^2|$. Since
$$(b_{1}-a_{1})^2 + (b_{2}-a_{2})^2 - (b_{1}-a_{2})^2 - (b_{2}-a_{1})^2 = 2(a_{2}-a_{1})(b_{1}-b_{2})$$
holds for every $a_{1},a_{2} \in A$, $b_{1},b_{2} \in B$, we have the inclusion
$$2(A-A)(B-B) + 2(A-A)(B-B) \subset 4\cdot (A-B)^2 - 4\cdot (A-B)^2.$$ 
On one hand Lemma 8 gives
\begin{eqnarray*}
|4(A-B)^2-4(A-B)^2| &\leq& \frac{|(A-B)^2+(A-B)^2|^{8}}{|(A-B)^2|^{7}} \\
&=& \frac{|(A-B)^2+(A-B)^2|^{8}}{|A-B|^{7}}.
\end{eqnarray*}
On the other hand, the above inclusion gives
\begin{eqnarray*}
|4 \cdot (A-B)^2-4 \cdot (A-B)^2|^{2} &\geq& |2(A-A)(B-B) + 2(A-A)(B-B) |^{2}\\
&=& |(A-A)(B-B)+(A-A)(B-B)|^{2}.
\end{eqnarray*}
Furthermore, Lemma 6 applied for $R=S=A-A$, $T=B-B$ tells us that
$$|(A-A)(B-B)+(A-A)(B-B)|^2 \gg \left|\frac{A-A}{A-A}\right| |B-B|^2.$$
By Lemma 7,
$$\left|\frac{A-A}{A-A}\right| \gg |A-A|^{1+\frac{1}{12}} \log^{-1/4}|A-A|,$$
so
$$|4 \cdot (A-B)^2-4 \cdot (A-B)^2|^{2} \gg |A-A|^{1+\frac{1}{12}} |B-B|^{2} \log^{-1/4}|A-A|.$$
We conclude that
$$\frac{|(A-B)^2+(A-B)^2|^{16}}{|A-B|^{14}} \gg |A-A|^{1+\frac{1}{12}} |B-B|^{2} \log^{-1/4}|A-A|.$$
By using Lemma 6 for $R=S=B-B$ and $T=A-A$ instead, we can similarly get
$$\frac{|(A-B)^2+(A-B)^2|^{16}}{|A-B|^{14}} \gg |A-A|^{2} |B-B|^{1+\frac{1}{12}} \log^{-1/4}|B-B|.$$
Putting everything together, we get
$$\min\left\{|A-A|, |B-B|, |A-B|\right\} \ll |(A-B)^2+(A-B)^2|^{1-\frac{13}{205}} \cdot \L(A,B),$$
where
$$\L(A,B)=\min\left\{ \log^{\frac{3}{205}} |A|, \log^{\frac{3}{205}} |B|\right\}.$$
This completes the proof. 
$\hfill \square$

\bigskip

\section{Some Remarks}

\bigskip

Theorem 2 is still far from being optimal. We conjecture that when $|\Delta(A \times A)| = o(|A|^{2})$,
$$|A-A| \ll |A|^{1+\epsilon},$$
for any $\epsilon > 0$. The $\epsilon$ in the conjecture is justfied by the following remark. Translate the set $A$ so that it contains $0$ so that now $|(A-A)^2+(A-A)^2| = o(|A|^{2})$ implies $|A^{2}+A^{2}| = o(|A|^2)$. On the other hand, using an argument similar to the one of Elekes and Ruzsa from [4], one can show that for every $A, B \subset \mathbb{R}$ we have that
$$|A^{2}+B^{2}| |A-A+B|^2 |A-A-B|^2 \gg |A|^{4} |B|^{2}.$$
In particular, whenever $A=B$ and $|A-A| \ll A$, Lemma 8 yields $|A^2+A^2| \gg |A|^{2}$, so we get a contradiction.

It is worth mentioning that even assuming the full-strength of the Erd\H{o}s-Szemer\'edi conjecture [6], which says that for any $\epsilon > 0$ one has
$$\max{|D^2+D^2|,|D^{2}D^{2}|} \gg |D|^{2-\epsilon'},$$
our proof for Theorem 2 only gives
$$|A-A| \ll |A|^{2-\frac{4}{7} + \epsilon}.$$
Using the updates building on the Konyagin and Shkredov improvements of Solymosi's bound [9,10], one can perhaps bring
$$|A-A| \ll |A|^{2-\frac{2}{7}}\log^{\frac{1}{7}} |A|$$
down to
$$|A-A| \ll |A|^{2-\frac{2}{7}-c}\log^{\frac{1}{7}} |A|$$
for some small constant $c>0$, but significant improvements to Theorem 2 should perhaps first come from replacing the inequality $|DD| \leq |2 \cdot D^{2} - 2 \cdot D^{2}|$ coming from $2DD \subset 2 \cdot D^{2} - 2 \cdot D^{2}$ with a more efficient argument. We believe that the Erd\H{o}s-Szemer\'edi conjecture should imply the claim that cartesian products with $o(|A|^{2})$ distinct distances satisfy $|A-A| \ll |A|^{1+\epsilon}$.

\bigskip

\subsection*{Acknowledgments} 

\bigskip

I would like to thank Oliver Roche-Newton, Misha Rudnev and Adam Sheffer for helpful conversations.

\bigskip

\small{ \textsc{California Institute of Technology, Pasadena, CA}

{\it{E-mail address}}: apohoata@caltech.edu}

\end{document}